\documentclass[10pt,twoside]{article}
\usepackage{graphicx}
\usepackage{amsmath}
\usepackage{Latex-document}
\usepackage{amssymb}

\newtheorem{theorem}{Theorem}[subsection]
\newtheorem{lemma}[theorem]{Lemma}
\newtheorem{corollary}[theorem]{Corollary}

\newtheorem{definition}[theorem]{Definition}

\def\Ker{{\mathrm {Ker}}}

\def\ban{{\mathrm {ban}}}
\def\Id{{\mathrm {Id}}}
\def\red{{\mathrm {red}}}

\def\R{{\mathbb R}}
\def\H{{\mathbb H}}
\def\Z{{\mathbb Z}}

\def\C{{\mathbb C}}
\def\K{{\mathcal K}}
\def\L{{\mathcal L}}
\def\A{{\mathcal A}}

\renewcommand{\thesection}{\arabic{section}.}
\renewcommand{\thesubsection}{\thesection\arabic{subsection}}

\def\addsec{\addtocounter{subsection}{1} \setcounter{theorem}{0}}

\title{\bf Banach {\boldmath $KK$}-theory and \vskip -2mm
the Baum-Connes Conjecture \vskip 6mm}
\author{V. Lafforgue\vspace*{-0.5cm}\thanks{Institut de Math\'ematiques
de Jussieu,  175 rue du Chevaleret, 75013 Paris, France. E-mail:
vlafforg@math.jussieu.fr}}

\date{\vspace{-8mm}}

\begin{document}
\maketitle

\thispagestyle{first} \setcounter{page}{795}

\begin{abstract}\vskip 3mm
The report below describes the applications of Banach KK-theory to a
conjecture of P. Baum and A. Connes about the K-theory of group
$C^*$-algebras, and a new proof of the classification by
Harish-Chandra, the construction by Parthasarathy
and the exhaustion by Atiyah and Schmid of the discrete series
representations of connected semi-simple Lie groups.

\vskip 4.5mm


\end{abstract}

\vskip 12mm

This report is intended to be very elementary. In the first part
we outline the main results in Banach KK-theory and the
applications to the Baum-Connes conjecture. In the second part we
show how the Baum-Connes conjecture for connected semi-simple Lie
groups can be applied to recover the classification of the
discrete series representations.

\section{Banach \boldmath{$KK$}-theory and the Baum-Connes conjecture}
\label{section 1}\setzero \vskip-5mm \hspace{5mm}

There are many surveys on Kasparov's $KK$-theory and the
Baum-Connes conjecture (see
\cite{baumconnes,SkExp,SkBki,higsonICM,julgBki,connes,valette})
and on
 Banach KK-theory
(\cite{SkBki,ecm}).

\subsection*{1.1. Generalized Fredholm modules} \addsec

\vskip-5mm \hspace{5mm}

We wish to define $A$-linear Fredholm operators (where $A$ is a
Banach algebra), with an index in $K_0(A)$. If $A=\C$, this index
should be the usual index of $\C$-linear Fredholm operators in
$K_0(\C)=\Z$.

We define a Banach algebra as a (non necessarily unital) $\mathbb
C$-algebra $A$ that is complete for a norm $\|.\|$ satisfying
$\|ab\|\leq \|a\| \|b\|$ for any $a,b \in A$.
If $A$ and $B$ are Banach algebras a morphism $\theta:A\to B$ is an algebra
morphism such that $\|\theta(a)\|\leq \|a\|$ for any $a\in A$.

$K_0$ and $K_1$ are two covariant functors from the category of Banach
algebras to the category of abelian groups.
 If $X$ is a locally compact space and $C_0(X)$ the algebra of
continuous functions vanishing at infinity, $K_0(C_0(X))$ and
$K_1(C_0(X))$ are the Atiyah-Hirzebruch K-theory groups.
For technical reasons we shall restrict ourselves to unital
Banach algebras in this subsection.

Let $A$ be  a unital Banach algebra.

A  right $A$-module $E$ is finitely generated projective  if
and only if
it is a direct summand in $A^n$ for some integer $n$. The set of
isomorphism classes of right finitely generated projective $A$-modules
is a semigroup because
the direct sum of two right finitely generated projective $A$-modules
is a right
finitely generated projective $A$-module. Then $K_0(A)$ is the
universal group associated to this semigroup ({\it i.e.} the group of formal
differences of elements of the semigroup). If $\theta:A\to B$ is a morphism
of unital Banach algebras, and $E$ is a right finitely
generated projective $A$-module then $E\otimes _A B$ is a right
finitely generated projective   $B$-module and this defines
 $\theta_*:K_0(A)\to K_0(B)$.

There is another definition of $K_0(A)$ for which the functoriality is
even more obvious~: $K_0(A)$ is the quotient of the free abelian group
generated by all idempotents $p$ in $M_k(A)$ for some integer $k$, by
the relations $\bigg[ \begin{pmatrix} p &
0\\0&q\end{pmatrix}\bigg]=[p]+[q]$
for any idempotents $p\in M_k(A)$ and $q\in M_l(A)$ and $[p]=[q]$ if
$p,q$ are idempotents of $M_k(A)$ and are connected by a path of
idempotents in $M_k(A)$ and $[0]=0$ where $0$ is the idempotent $0$ in
$M_k(A)$. The link with the former definition is that any
idempotent $p\in M_k(A)$ acts on the left on $A^k$ as a projector $P$
and $\mathrm{Im} P$ is a right finitely generated projective
$A$-module (it is a direct summand in the right $A$-module $A^k$).

The following construction was performed for $C^*$-algebras by
Mischenko and Kasparov, in connection with the Novikov conjecture
(\cite{mishchenko2,kaspnov}). We adapt it to Banach algebras.

%
A right Banach $A$-module is a Banach space (with a given norm $\|.\|_E$)
equipped with a right action of $A$ such that $1\in A$ acts by
identity and $\|xa\|_E\leq \|x\|_E\|a\|_A$ for any $x\in E$ and $a\in
A$.  Let $E$ and $F$ be right Banach $A$-modules. A morphism $u:E\to
F$ of right Banach
$A$-modules is a continuous $\C$-linear map such that $u(xa)=u(x)a$ for any
$x\in E$ and $a\in A$. The space $\mathcal L_A(E,F)$ of such morphisms
is a Banach space with norm $\|u\|=\sup_{x\in E,
\|x\|_E=1}\|u(x)\|_F$. A morphism $u\in \L_A(E,F)$ is said to be ``$A$-rank
one'' if $u=w\circ v$ with $v\in \L_A(E,A)$ and $w\in \L_A(A,F)$. The
space $\mathcal K_A(E,F)$
of $A$-compact morphisms is the closed vector span of $A$-rank one morphisms
in $\L_A(E,F)$. If $E=F$, $\L_A(E)=\L_A(E,E)$ is a Banach algebra and
$K_A(E)=K_A(E,E)$ is a closed ideal in it.

\begin{definition} A Fredholm module over $A$ is the data of a $\mathbb
Z/2$ graded  right Banach $A$-module  $E$ and an odd morphism $T\in
\mathcal L_A(E)$ such that $T^2-\mathrm{Id}_E\in \K_A(E)$.
\end{definition}

In other words $E=E_0\oplus E_1$, $T=\begin{pmatrix} 0 & v\\u
&0\end{pmatrix}$ and $u\in \mathcal L_A(E_0,E_1)$ and $v\in \mathcal
L_A(E_1,E_0)$ satisfy $vu-\mathrm{Id}_{E_0}\in \mathcal K_A(E_0)$ and
$uv-\mathrm{Id}_{E_1}\in \mathcal K_A(E_1)$.

If $(E,T)$ is a Fredholm module over $A$ and $\theta:A\to B$ a
unital morphism then $(E\otimes _A B,T\otimes 1)$ is a
Fredholm module over $B$ (here $E\otimes _A B$ is the completion of
$E\otimes _A^{\mathrm{alg}} B $ for the maximal Banach norm such that
$\|x\otimes b\|\leq \|x\|_E \|b\|_B$ for $x\in
E$ and $b\in B$).

Let $A[0,1]$ be the Banach algebra of continuous functions from $[0,1]$
to $A$ with the norm $\|f\|=\sup_{t\in [0,1]}\|f(t)\|_A$ and
$\theta_0,\theta_1:A[0,1]\to A$ the
evaluations at $0$ and $1$. Two Fredholm modules on $A$ are
said to be homotopic if they are the images by $\theta_0$ and $\theta_1$
of a Fredholm module over $A[0,1]$.

\begin{theorem}\label{fredholm} There is a functorial bijection
between $K_0(A)$ and the
set of homotopy classes of Fredholm modules over $A$, for any unital
Banach algebra $A$.
\end{theorem}

Let $(E_0,E_1,u,v)$ be a Fredholm module over $A$. Its index,
{\it i.e.} the corresponding element in $K_0(A)$, is constructed
as follows. It is possible to find $n\in \mathbb N$ and $w\in
\mathcal K_A( A^n,E_1)$ such that  $(u,w)\in \mathcal
L_A(E_0\oplus A^n,E_1)$  is surjective. Its kernel is then
finitely generated projective  and the index is the formal
difference  of $\mathrm{Ker}((u,w))$ and $A^n$.

An ungraded Fredholm module over $A$ is the data of a (ungraded) right
Banach module $E$ over $A$, and $T\in \L_A(E)$ such that
$T^2-\Id_E\in \K_A(E)$. There is a functorial bijection between
$K_1(A)$ and the set of homotopy classes of ungraded Fredholm modules.

For a non-unital algebra $A$,
$K_0(A)=\Ker(K_0(\tilde A)\to K_0(\C)=\Z)$ and $K_1(A)=K_1(\tilde A)$
where $\tilde A=A\oplus \C 1$. In particular every idempotent in
$M_k(A)$ gives a class in $K_0(A)$ but in general not all classes in
$K_0(A)$ are obtained in this way.
The definition of a Fredholm module should be slightly modified for
non-unital Banach algebras, but the theorem~\ref{fredholm} remains
true.

\subsection*{1.2. Statement of the Baum-Connes conjecture} \addsec

\vskip-5mm \hspace{5mm}

Let $G$ be  a second countable, locally compact group. We fix a
left-invariant Haar measure $dg$ on $G$.
Denote by $C_c(G)$ the convolution algebra of complex-valued
continuous compactly
supported functions on $G$. The convolution of $f,f' \in C_c(G)$ is
given by
$f*f'(g)=\int _G f(h)f'(h^{-1}g)dh$ for any $g\in G$.

When $G$ is discrete and $dg$ is the counting measure,
$C_c(G)$ is also denoted by $\C G$ and if $e_g$ denotes the delta
function at $g\in G$ (equal to $1$ at $g$ and $0$ elsewhere),
$(e_g)_{g\in G}$ is a basis of $\C G$ and the convolution product is
given by $e_ge_{g'}=e_{gg'}$.

The completion of $C_c(G)$ for the norm $\|f\|_{L^1}=\int _G|f(g)|dg$
is a Banach algebra and is denoted by $L^1(G)$.

For any $f\in C_c(G)$ let $\lambda(f)$ be the operator $f'\mapsto
f*f'$ on $L^2(G)$.
The completion of $C_c(G)$ by the operator norm
$\|f\|_{\red}=\|\lambda(f)\|_{\L_{\C}(L^2(G))}$ is called the reduced
$C^*$-algebra of $G$ and denoted by $C^*_{\red}(G)$. If $G$ is
discrete $(e_{g'})_{g'\in G}$ is an orthonormal basis of $L^2(G)$ and
$\lambda (e_g):e_{g'}\mapsto e_{gg'}$.

%

For any $f\in C_c(G)$, $\|f\|_{L^1}\geq \|f\|_{\mathrm{red}}$ and
$L^1(G)$ is a dense subalgebra of $C^*_{\mathrm{red}}(G)$. We denote
by $i:L^1(G)\to C^*_{\mathrm{red}}(G)$ the inclusion.
%
%
\vskip 0.5cm

Assume now that $M$ is a smooth compact manifold, and $\tilde M$ a
Galois covering of $M$ with group $G$ (if $\tilde M$ is simply
connected, $G=\pi_1(M)$). Let $E_0$ and $E_1$ be two smooth hermitian
finite-dimensional vector bundles over $M$
and $u$ an order $0$
elliptic pseudo-differential operator from $L^2(M,E_0)$ to
$L^2(M,E_1)$. Since $u$ is elliptic there is an order $0$
pseudo-differential operator  $v:L^2(M,E_1)\to L^2(M,E_0)$ such that
$\mathrm{Id}_{L^2(M,E_0)}-vu$ and $\mathrm{Id}_{L^2(M,E_1)}-uv$ have
order $\leq -1$ and therefore are compact.
Let $\mathcal E$ be the quotient of $\tilde M \times C^*_{\mathrm{red}}(G)$ by the
diagonal action of $G$ ($G$ acting on $C^*_\red(G)$ by left
translations)~: $\mathcal E$ is a flat bundle of right
$C^*_\red(G)$-modules over $M$, whose fibers are isomorphic to
$C^*_{\mathrm{red}}(G)$. Then $L^2(M,E_0\otimes \mathcal
E)$ and $L^2(M,E_1\otimes \mathcal
E)$ are right Banach (in fact Hilbert) modules over $C^*_{\mathrm{red}}(G)$
and it is possible to
lift $u$ and $v$ to $\tilde u$ and $\tilde v$ so that
$(L^2(M,E_0\otimes \mathcal E),L^2(M,E_1\otimes \mathcal E), \tilde
u,\tilde v)$ is a Fredholm module over $C^*_{\mathrm{red}}(G)$, whose
index lies in $K_0(C^*_{\mathrm{red}}(G))$ and the index does not
depend on the choice of the liftings.


The operator $u$ represents a ``$K$-homology class'' in $K_0(M)$, and
using the classifying map $M\to BG$, it defines  an element of $K_{0,c}(BG)$,
the K-homology with compact support of the classifying space $BG$. For
any discrete
group $G$ we can define a morphism of abelian groups  $K_{*,c}(BG)\to
K_*(C^*_{\mathrm{red}}(G))$ ($*=0,1$).
This morphism is the Baum-Connes assembly map when $G$ is discrete and
torsion free. When $G$ is not discrete or has
torsion, the index construction can be performed starting from
a proper action of $G$ (instead of the free and proper action of $G$ on $\tilde
M$ in the last paragraph), and therefore we have to introduce the
 space $\underline{E}G$ that classifies the proper actions of $G$.
Using Kasparov equivariant KK-theory, the $G$-equivariant K-homology
$K_*^G(\underline{E}G)$ with $G$-compact support ($*=0,1$) may be
defined, and there is an assembly map
$$\mu_{\mathrm{red}}: K_*^G(\underline{E}G)\to K_*(C^*_{\mathrm{red}}(G)). $$
In the same way we can define
$\mu_{L^1} : K_*^G(\underline{E}G)\to K_*(L^1(G)) $
and  $\mu_{\mathrm{red}}=i_*\circ \mu_{L^1}$.
\vskip 0.4cm
\noindent
{\bf Baum-Connes conjecture }\cite{preprintBaumConnes,baumconnes} :
 {\it If $G$ is a second countable, locally
compact group then the assembly map $\mu_{\mathrm{red}}:
K_*^G(\underline{E}G)\to  K_*(C^*_{\mathrm{red}}(G))$ is an
isomorphism.}
 \vskip 0.4cm
\vskip 0.2cm

 Bost conjectured~: If $G$ is a second countable, locally
compact group (and has reasonable geometric properties)
then the assembly map $\mu_{L^1} :
K_*^G(\underline{E}G)\to K_*(L^1(G))$ is an isomorphism.
\vskip 0.2cm

In many cases  $K_*^G(\underline{E}G)$ can be computed. For instance
if $G$ is a
discrete torsion free subgroup of a reductive Lie group $H$ and $K$ is
a maximal compact subgroup of $H$, then a possible $\underline{E}G$ is
 $ H/K$  and $K_*^G(\underline{E}G)$ is the K-homology with compact support
of $G\backslash H/K$. This group may be computed
thanks to Mayer Vietoris sequences. See part 2 for the case where
$G$ is a Lie group.

\subsection*{1.3. {\boldmath $KK$}-theory} \addsec

\vskip-5mm \hspace{5mm}

For any $C^*$-algebras $A$ and $B$, Kasparov~\cite{kaspnov,kaspkk}
defined  an abelian group $KK(A,B)$, covariant in $B$ and
contravariant in $A$. There is a  product
$KK(A,B)\otimes KK(B,C)\to KK(A,C)$. Moreover $KK(\C,A)=K_0(A)$ and
therefore the product gives a morphism $KK(A,B)\to
\mathrm{Hom}(K_0(A),K_0(B))$. The definition of $KK(A,B)$ is like
definition~\ref{defikk} below, but with Hilbert modules instead of
Banach modules.

For any Banach algebras $A$ and $B$, we define~\cite{these,SkBki} an abelian group \linebreak $KK^\ban(A,B)$,
covariant in $B$ and contravariant in $A$. There is no product, but a morphism $KK^\ban(A,B)\to
\mathrm{Hom}(K_0(A),K_0(B))$. Assume that $B$ is unital (otherwise the definition has to be slightly modified).

\begin{definition}\label{defikk}
$E^\ban(A,B)$ is the set of isomorphism classes of  data
$(E,\pi,T)$, where $E$ is a
$\Z/2\Z$-graded right Banach module, $\pi: A\to \L_B(E)$ is a morphism
of Banach algebras and takes values in even operators, and $T\in
\L_B(E)$ is odd and satisfies $a(T^2-\Id_E)\in K_B(E)$ and $aT-Ta\in
K_B(E)$ for any $a\in A$.
\end{definition}

Then $KK^\ban(A,B)$ is the set of homotopy classes in  $E^\ban(A,B)$,
where the homotopy
relation is defined using $E^\ban(A,B[0,1])$.

\noindent {\bf Remark :}
$E^\ban(\C,B)$ is the set of isomorphism classes of Fredholm
modules over $B$ and $KK^\ban(\C,B)=K_0(B)$.

If $p$ is an idempotent in $A$, and $(E,\pi,T)\in E^\ban(A,B)$, the
image of $[p]\in K_0(A)$ by the image of $[(E,\pi,T)]\in KK^\ban(A,B)$
in $\mathrm{Hom}(K_0(A),K_0(B))$ is defined to be the index of the
Fredholm module over $B$ equal to $(\mathrm{Im}
\pi(p),\pi(p)T\pi(p))$. When $p$ is an idempotent in $M_k(A)$, we use
the image of $p$ by $M_k(A)\to \L_B(E^k)$. This is enough to define
the morphism $KK^\ban(A,B)\to \mathrm{Hom}(K_0(A),K_0(B))$, when $A$
is unital.

The same definition with ungraded modules gives $KK_1^\ban(A,B)$, and, with
the notation $KK=KK_0$, we have a morphism
$KK_i^\ban(A,B)\to \mathrm{Hom}(K_j(A),K_{i+j}(B))$, where all the indices
are modulo $2$.

\subsection*{1.4. Status of injectivity and the element {\boldmath $\gamma$}} \addsec

\vskip-5mm \hspace{5mm}

The injectivity of the Baum-Connes map $\mu_{\mathrm{red}}$ (and therefore of
$\mu_{L^1}$) is known for the following very large classes of groups~:
\begin{trivlist}
\item  a) groups acting continuously properly isometrically on a
complete simply
connected riemannian manifold with controlled
non-positive sectional curvature,
and in particular  closed subgroups of reductive Lie groups
(\cite{conspectus,kaspkk}),
\item  b) groups acting continuously properly isometrically on an
affine building and in particular  closed subgroups of reductive
$p$-adic groups (\cite{immeublesnovikov}),
\item  c) groups acting continuously properly isometrically on a
discrete metric space with good properies at infinity (weakly geodesic,
uniformly locally finite, and ``bolic'' \cite{crasskandkasp,bolic}), and in
particular hyperbolic groups ({\it i.e.} word-hyperbolic in the sense
of Gromov),
\item  d)  groups acting continuously amenably on a compact space
(\cite{Higson}).
\end{trivlist}

In the cases a),b),c) above, the proof of injectivity provides an
explicit idempotent endomorphism on $K_*(C^*_{\mathrm{red}}(G))$ whose
image is the image of $\mu_{\mathrm{red}}$ (and the same for $\mu_{L^1}$).
In case d), J.-L. Tu has also constructed such an endomorphism, but in
a less explicit way.

To state this we need to understand a baby case of Kasparov's
equivariant KK-groups. Let $G$ be a second countable, locally compact group.
We denote by  $E_G(\C,\C)$ the set of isomorphism classes of triples
$(H,\pi,T)$ where $H$ is a $\mathbb Z/2$-graded Hilbert space, $\pi$ a
unitary representation of $G$ on $H$ (such that for any $x\in H$,
$g\mapsto gx$ is continuous from $G$ to $H$) and $T$ an odd operator
on $H$ such
$T^2-\mathrm{Id}_H$ is compact and $\pi(g)T\pi(g^{-1})-T$ is compact
and depends norm continuously on $g\in G$.
Then $KK_G(\mathbb C,\mathbb C)$ is the quotient of $E_G(\C,\C)$ by homotopy.
Kasparov proved that $KK_G(\mathbb C,\mathbb C)$ has a ring structure
(using direct sum for the addition and tensor products together with a
quite difficult construction for the multiplication).

If $\pi$ is a unitary representation of $G$ on a Hilbert space $H_0$
and $H_1=0$ then  $(H,\pi,0)\in E_G(\C,\C)$ if and only if $H_0$ has
finite dimension. If moreover
$H_0=\C$ and $\pi$ is the trivial  representation of $G$, the
class of $(H,\pi,0)$ is the unit of
$KK_G(\C,\C)$ and is denoted by $1$. If $G$ is compact the classes of
 $(H,\pi,0)$  with $H_1=0$ (and $\dim
H_0<+\infty$) generate $KK_G(\mathbb{C},\mathbb{C})$ and
$KK_G(\mathbb{C},\mathbb{C})$ is equal to the representation ring of $G$.

The important fact is that there is a ``descent morphism''
$$j_{\mathrm{red}}: KK_G(\mathbb C,\mathbb C)\to \mathrm{End}(K_*(C^*_{\mathrm{red}}(G))). $$
In fact it is a ring homomorphism and
 $j_{\mathrm{red}}(1)=\mathrm{Id}_{K_*(C^*_{\mathrm{red}}(G))}$.
It is defined as the composite of two maps
 $KK_G(\mathbb C,\mathbb C)\to KK(C^*_{\mathrm{red}}(G),C^*_{\mathrm{red}}(G)) \to
\mathrm{End}(K_*(C^*_{\mathrm{red}}(G)))$. The construction of
 $j_{\mathrm{red}}$ is due to Kasparov. The construction of $j_{L^1}$
to be explained below is an adaptation of it.

The following extremely important theorem also contains earlier works
of Mishchenko and Solovjev.
\begin{theorem}\label{gammainj}(Kasparov, Kasparov-Skandalis
\cite{kaspkk,immeublesnovikov,crasskandkasp,bolic})
If $G$ belongs to one of the classes a),b),c) above, the geometric
conditions in a),b) or c) allow
to construct an idempotent element $\gamma\in KK_G(\mathbb C,\mathbb
C)$ such that $\mu_{\mathrm{red}}$
is injective and its image is equal to the image of the idempotent
 $j_{\mathrm{red}}(\gamma)\in
\mathrm{End}(K_*(C^*_{\mathrm{red}}(G)))$.
\end{theorem}

\subsection*{1.5. Homotopies between {\boldmath $\gamma$} and {\boldmath $1$}} \addsec

\vskip-5mm \hspace{5mm}

We assume that $G$ belongs to one of the classes a),b),c). Then
the injectivity of $\mu_{\mathrm{red}}$ is known  and the surjectivity
is equivalent to the equality $j_{\mathrm{red}}(\gamma)=\mathrm{Id}
\in
\mathrm{End}(K_*(C^*_{\mathrm{red}}(G)))$.

\begin{theorem} We have $\gamma=1$ in $KK_G(\mathbb C,\mathbb C)$ if
\begin{enumerate}
\item $G$ is a  free group (Cuntz, \cite{cuntzKamenability}) or a
closed subgroup of $SO(n,1)$ (Kasparov, \cite{lorentz})
or of $SU(n,1)$ (Julg-Kasparov, \cite{JulgKasp}) or of $SL_2(\mathbb F)$
with $\mathbb F$ a local non-archimedian field (Julg-Valette,
\cite{julgvalette}),
 \item $G$ acts isometrically and properly on a Hilbert space
(Higson-Kasparov \cite{HigsonKasp,julgBki}).
\end{enumerate}
\end{theorem}

In fact the second case contains the first one.

If $G$ has property (T) and is not
compact, $\gamma \neq 1$ in $KK_G(\C,\C)$ :  it is impossible to
deform $1$ to $\gamma$ in $E_G(\C,\C)$ because
the trivial representation is isolated among unitary representations
of $G$ if $G$ has
property (T) and  $\gamma$ can be represented by $(H,\pi,T)$ such that
$H$ has no invariant vector (and even $H$ is tempered). All simple
real or $p$-adic groups of rank $\geq 2$, and $Sp(n,1)$ and
$F_{4 (-20)}$, and all their lattices,  have property (T)
(see~\cite{delaharpevalette}).

It is then natural to broaden the class of representations  in order
to break the isolation of the
trivial one. In \cite{julgrepunifbornees} Julg proposed to use
uniformly bounded representations on Hilbert spaces (to solve the case
of $Sp(n,1)$).

For any non compact group $G$ the trivial representation is not
isolated among isometric representations in Banach spaces (think of
the left regular representation on $L^p(G)$, $p$ going to infinity).

\begin{definition} Let $E_G^\ban(\C,\C)$ be the set of isomorphism
classes of  triples \linebreak $(E,\pi,T)$ with $E$ a $\mathbb Z/2$-graded Banach space endowed with an isometric
representation of $G$ (such that $g\mapsto gx$ is continuous from $G$ to $E$ for any $x\in E$),  $T\in \mathcal
L_{\C}(E)$ an odd operator such that $T^2-\mathrm{Id}_E$ belongs to $\mathcal K_{\mathbb C}(E)$ and
$\pi(g)T\pi(g^{-1})-T$ belongs to $\mathcal K_{\mathbb C}(E)$ and depends norm continuously on $g\in G$.
\end{definition}

Then $KK^{\mathrm{ban}}_G(\mathbb C,\mathbb C)$ is defined as the
quotient of $E_G^\ban(\C,\C)$ by  homotopy.

Since any unitary representation of $G$ on a Hilbert space $H$ is an
isometric representation on the Banach space $H$,
there is a natural morphism of abelian groups
$KK_G(\mathbb C,\mathbb C)\to KK^{\mathrm{ban}}_G(\mathbb C,
\mathbb C)$.

To state our main theorem, we need to look at slightly smaller classes
than a) and c) above. We call these new classes a') and c'). They are
morally the same, and in particular they respectively contain all
closed subgroups of reductive Lie groups, and all hyperbolic groups
(for general hyperbolic groups see~\cite{mineyevyu}, and~\cite{these}
for a slightly different approach).

\begin{theorem}\cite{these,SkBki}
For any group $G$ in the classes a'),  b), or c'), we have $\gamma=1$ in
$KK^{\mathrm{ban}}_G(\mathbb C, \mathbb C)$.
\end{theorem}

In fact the statement is slightly incorrect, we should allow
representations with a slow growth, but this adds no real
difficulty.
The proof of this theorem is quite technical. Let me just indicate
some ingredients involved.
If $G$ is in class a') then $G$ acts continuously isometrically properly on a
complete simply connected riemmannian manifold $X$ with controlled non-positive
sectional curvature, and $X$ is contractible (through geodesics) and
the de Rham cohomology of $X$ (without support) is $\mathbb C$ in
degree $0$ and $0$ in other degrees. It is possible to put norms on
the spaces of differential forms (on which $G$ acts) and to build a
parametrix for
the de Rham operator (in the spirit of the Poincar\'e lemma) in order to
obtain a resolution of the trivial
representation, and in our language an element of $E_G^\ban(\C,\C)$
 equal to $1$ in $KK^{\mathrm{ban}}_G(\mathbb C, \mathbb
C)$. The norms we use are essentially Sobolev $L^{\infty}$ norms. Then it
is possible to conjugate the operators by an exponential of the
distance to a fixed point in $X$ and then to deform these norms to Hilbert norms
(through $L^p$ norms, $p\in [2,+\infty]$) and to reach $\gamma$.

If $G$ belongs to class b) the de Rham complex is replaced by the
simplicial homology complex (with $L^1$ norms) on the
building. If $G$ belongs to class c') a Rips complex plays the same
role as the building in b).

It is not possible to apply directly this theorem to the Baum-Connes
conjecture
because there is no obvious descent map
$KK^{\mathrm{ban}}_G(\mathbb C, \mathbb C)\to
\mathrm{End}(K_*(C^*_{\mathrm{red}}(G)))$, and in the next subsection we
shall see the difficulties encountered and the way one bypasses them
in a few cases.

On the other hand, we  may apply this theorem to Bost conjecture,
because there is descent map $j_{L^1}:KK_G^\ban(\C,\C)\to
KK^\ban(L^1(G),L^1(G))$.

We explain it when $G$ is discrete. Let $(E,\pi,T)\in
E_G^\ban(\C,\C)$. We denote by $L^1(G,E)$ the
completion of
 $E\otimes \C G$ for the norm $\|\sum_{g\in G}x(g)\otimes e_g\|=\sum_{g\in G}
\|x(g)\|_E$. Then $L^1(G,E)$ is a right Banach $L^1(G)$-module by the formula  $(x\otimes e_g)e_{g'}=x\otimes
e_{gg'}$ and there is a Banach algebra morphism $\hat \pi:L^1(G)\to \L_{L^1(G)}(L^1(G,E))$ by the formula $\hat
\pi(e_{g'})(x\otimes e_g)=\pi(g')(x)\otimes e_{g'g}$. Then $(L^1(G,E), \hat \pi$, $T\otimes 1)\in
E^\ban(L^1(G),L^1(G))$ gives the desired class in $KK^\ban(L^1(G),L^1(G))$.

This and section 1.3 imply   the Bost conjecture in many cases.

\begin{theorem}\label{bostconj}
For any group $G$ in the classes a'), b) or c'),
$\mu_{L^1}:K_*^G(\underline{E}G)\to K_*(L^1(G))$ is an
isomorphism.
\end{theorem}

\subsection*{1.6. Unconditional completions}\label{unconditional} \addsec

\vskip-5mm \hspace{5mm}

Let $G$ be a second countable, locally compact group. Let
$\mathcal A(G)$ be a Banach algebra containing
$C_c(G)$ as a dense subalgebra. We write
$\mathcal A(G)$ instead of $\mathcal A$ for notational convenience.
We ask for a necessary and sufficient condition such that there is a
``natural'' descent map $j_{\mathcal A}:KK^{\mathrm{ban}}_G(\mathbb
C, \mathbb C)\to KK^\ban(\mathcal A(G),\mathcal A(G))$.

In order to simplify the argument below, we will assume $G$ to be
discrete.

Let $E$ be a Banach space with an isometric representation of $G$.
Then $E\otimes \C G$ has a right $\C G$-module structure given by
$(x\otimes e_g)e_{g'}=x\otimes e_{gg'}$ and there is a morphism
$\hat \pi :\C G\to \mathrm{End}_{\C G}(E\otimes \C G)$ given by the formula
  $\hat \pi(e_{g'})(x\otimes e_g)=\pi(g')(x)\otimes e_{g'g}$.
We look for a completion $\A(G,E)$ of $E\otimes \C G$ by a Banach norm
such that $\A(G,E)$ is a right Banach $\A(G)$-module
and $\hat \pi$ extends to a morphism of Banach algebras
$\hat \pi: \A(G) \to \L_{\A(G)}(\A(G,E))$.

In order to have enough $\A(G)$-rank one operators, it is quite
natural to assume that the norm on $\A(G,E)$ satisfies~: for any $x\in
E$ and $\xi\in
\mathcal L_{\C}(E,\mathbb C)$, if we denote by
$R_{x}:\C G \to E\otimes \C G$ the map
$e_g \mapsto x \otimes e_g$ and by $S_{\xi}:E\otimes \C G \to \C G$
the map $y \otimes e_g \mapsto \xi(y)e_g$, we have
$\|R_x(f)\|_{\A(G,E)}\leq \|x\|_{E}\|f\|_{\A(G)}$
for any $f\in \C G$ and
$\|S_{\xi}(\omega)\|_{\A(G)}\leq
\|\xi\|_{\L_{\C}(E,\C)}\|\omega\|_{\A(G,E)}$
for any $\omega\in E\otimes \C G$.
Now fix $x\in E$ and $\xi \in \L_{\C}(E,\mathbb C)$ and denote by $1$ the
unit in $G$. For any
$f=\sum_{g\in G}f(g)e_g \in
\C G$, $S_{\xi}(\hat \pi(f)(R_x( e_1)))$ is
$\sum_{g\in G} \xi(\pi(g)(x))f(g)e_g$ in $\C G$.
 For any function $c$ on $G$, we define  the Schur multiplication by
$c$ to be  the pointwise product $\C G\to \C G$, $\sum _{g\in
G}f(g)e_g\mapsto \sum_{g\in    G}c(g)f(g)e_g$.

In this way we obtain the following necessary condition~: for any $x\in E$ and
$\xi\in \mathcal L_{\C}(E,\mathbb C)$ the Schur multiplication by the matrix
coefficient $g\mapsto \xi(\pi(g)(x))$ is bounded from $\mathcal A(G)$
to itself
and its norm (in $\L_\C(\A(G))$) is less than
$\|x\|_E\|\xi\|_{\L_{\C}(E,\C)}$. But for any $L^\infty$-function
$c$ on $G$ we can find an isometric representation $\pi$ of $G$ on a
Banach space $E$ and $x\in E$ and $\xi\in \L_{\C}(E,\C)$ such that
$\|x\|_E\|\xi\|_{\L_{\C}(E,\C)}= \|c\|_{L^\infty}$  and
$c(g)=\xi(\pi(g)x)$
for any $g\in G$ (take $E=L^1(G)$, $x=\delta _1$, $\xi=c$).  Therefore
a necessary condition is that $\mathcal
A(G)$ is an unconditional completion in the following sense.

\begin{definition} A Banach algebra $\mathcal A(G)$ (with a given norm
$\|.\|_{\mathcal A(G)}$)
containing $C_c(G)$ as a dense subalgebra is called an unconditional
completion if the norm $\|f\|_{\A(G)}$ of $f\in C_c(G)$ only depends on
$g\mapsto |f(g)|,\  G\to \R_+$. \end{definition}

Remark that $L^1(G)$ is  an unconditional completion of $C_c(G)$
but  $C^*_{\red}(G)$ is not.

In fact this condition is also sufficient to construct the descent
map. For the sake of simplicity, we still assume that $G$ is
discrete. If $\mathcal A(G)$ is an
unconditional completion of $\C G$, and $(E,\pi,T)$ is in
$E_G^{\ban}(\C,\C)$, we define $\mathcal A(G,E)$ as the
completion of $E\otimes \C G$ for the norm $\|\sum _{g\in G}x(g)\otimes
e_g\|=\|\sum _{g\in G}\|x(g)\|_E \ e_g\|_{\mathcal A(G)}$ and $\mathcal
A(G,E)$ is a right Banach module over $\mathcal A(G)$ and there is a
morphism $\hat \pi: \A(G) \to \L_{\A(G)}(\A(G,E))$, and $(\A(G,E),
\hat \pi, T\otimes 1)\in E^\ban(\A(G),\A(G))$.

In this way, for any unconditional completion $\A(G)$ of $\C G$, we
have a descent map $j^{\ban}_{\A}:KK_G^{\ban}(\C,\C)\to
KK^\ban(\A(G),\A(G))\to \mathrm{End}(K_*(\A(G)))$. We can also define
an assembly map  $\mu_{\A}:K_*^G(\underline{E}G)\to K_*(\A(G))$. If
$\A(G)$ is an involutive subalgebra of $C^*_{\mathrm{red}}(G)$, and
$i:\A(G)\to C^*_\red(G)$ denotes the inclusion, $\mu_\red=i_*\circ \mu_\A$.

\begin{theorem}\label{uncondthm} (\cite{these})
For any group $G$ in the classes a'), b) or c'),
and for any unconditional completion $\A(G)$ of $C_c(G)$,
$\mu_{\A}:K_*^G(\underline{E}G)\to K_*(\A(G))$ is an isomorphism.
\end{theorem}

Let $A,B$ be Banach algebras and $i:A\to B$ an injective morphism of
Banach algebras. We say that $A$ is stable under holomorphic
functional calculus in $B$ if any element of $A$ has the same spectrum
in $A$ and in $B$. If $A$ is dense and stable
under holomorphic functional calculus in $B$ then $i_*:K_*(A)\to
K_*(B)$ is an isomorphism (see the appendix of \cite{bost}).

\begin{corollary}\label{maincor}
For any group $G$ in the classes a'),  b) or c'),
if $C_c(G)$ admits an unconditional completion $\A(G)$ which is  an
involutive subalgebra of $C^*_{\red}(G)$ and is stable under
holomorphic functional calculus in $C^*_{\red}(G)$,
then $\mu_{\red}:K_*^G(\underline{E}G)\to K_*(C^*_{\red}(G))$ is an
isomorphism.

This condition is fulfilled for
\begin{trivlist}
\item a) hyperbolic groups,
\item b) cocompact lattices in a product of a finite number of
groups among Lie or
$p$-adic groups of rank one,  $SL_3(\mathbb F)$ with $\mathbb F$ a
 local field (even $\H$) and  $E_{6(-26)}$,
\item c) reductive Lie groups and reductive groups over  non-archimedian
local fields.
\end{trivlist}
\end{corollary}

In case c), $\A(G)$ is a variant of the Schwartz algebra of the group
(\cite{these}). In this case the Baum-Connes conjecture was already
known for linear connected reductive groups (Wassermann
\cite{wassermann}) and for the $p$-adic $GL_n$ (Baum, Higson, Plymen
\cite{glpadique}).
In case a),b) this result is based on a property first introduced
by Haagerup for the free group and called (RD) (for rapid decay) by
Jolissaint (\cite{jolissaint}).
In case a),b) G has property (RD)~: this is due to
Haagerup for free groups (\cite{haagerup}), Jolissant for ``geometric
hyperbolic groups'', de la Harpe for general hyperbolic groups
(\cite{delaharpe}), Ramagge, Robertson and Steger for $SL_3$ of a
non-archimedian local field
(\cite{robertson}),   the author for $SL_3(\R)$ and $SL_3(\C)$
(\cite{sl3}), Chatterji for $SL_3(\H)$ and
$E_{6(-26)}$ (\cite{chatterji}), and the
remark that it holds for products is due to Ramagge, Robertson and
Steger (\cite{robertson}) in a particular case, and independantly to
Chatterji (\cite{chatterji}) and Talbi (\cite{talbi})
in general. A discrete
group $G$ has property (RD) if there is a lenght function $\ell :G\to \R_+$
({\it i.e.} a function satisfying $\ell (g^{-1})=\ell (g)$ and
$\ell (gh)\leq \ell (g)+\ell(h)$ for any $g,h\in G$) such that for
$s\in \R_+$ big enough, the completion $H^s(G)$ of $\C G$ for the norm
$\|\sum f(g)e_g\|_{H^s(G)}=\|\sum (1+\ell(g))^sf(g)e_g\|_{L^2(G)}$ is
contained
in $C^*_{\red}(G)$. Then, for $s$ big enough, $H^s(G)$
 is a Banach algebra and an involutive   subalgebra of
$C^*_{\red}(G)$ and is dense and
stable under holomorphic functional calculus (\cite{jolissaint,sl3});
 it is obvious that $H^s(G)$ is an unconditional completion of $\C
G$.

As a consequence of this result the Baum-Connes conjecture has been proven
for all almost connected groups by Chabert, Echterhoff and Nest
(\cite{chabertechtnest}).

\subsection*{1.7. Trying to push the method further} \addsec

\vskip-5mm \hspace{5mm}

In order to prove new cases of the surjectivity of the Baum-Connes map
(when the injectivity is proven and the $\gamma$ element exists) we
should look for a
dense subalgebra $\A(G)$ of $C^*_{\red}(G)$ that is stable under
holomorphic functional calculus and a homotopy between $\gamma$ and
$1$ through (perhaps special kind of) elements of $E^\ban_G(\C,\C)$
which all give a map $K_*(\A(G))\to K_*(C^*_{\red}(G))$ by
the descent construction. Thanks to the discussion in
subsection~\ref{unconditional} a necessary condition for this is that for
any $(E,\pi,T)$ in the homotopy between $\gamma$ and
$1$,  for any $x\in E$ and $\xi \in
\L_{\C}(E,\C)$, the Schur multiplication by the matrix coefficient
$g\mapsto \xi(\pi(g)(x))$ is bounded from $\A(G)$ to $C^*_{\red}(G)$
and has
norm $\leq \|x\|_E\|\xi\|_{\L_{\C}(E,\C)}$. So we should first look
for a homotopy between $\gamma$ and
$1$ such that the fewest possible  matrix coefficients
appear. For groups acting properly on buildings,
 this homotopy can be shown to exist. The problem
for general discrete groups properly acting on buildings is to find a
subalgebra $\A(G)$ of $C^*_{\red}(G)$ that is stable under
holomorphic functional calculus and satisfies the condition with
respect to these matrix coefficients. The first step (the crucial one
I think) should be to find
a subalgebra $\A(G)$ of $C^*_{\red}(G)$ that is stable under
holomorphic functional calculus and satisfies the following condition
: there is a integer $n$, a distance $d$ on the building and a point
$x_0$ on the building such that the Schur product by the
characteristic function of $\{ g\in G, d(x_0,gx_0)\leq r\}$ from $\A(G)$ to
$C^*_{\red}(G)$ has norm less than $(1+r)^n$, for any $r\in \R_+$.

\subsection*{1.8. The Baum-Connes conjecture with coefficients} \addsec

\vskip-5mm \hspace{5mm}

Let $G$ be a second countable, locally compact group and $A$ a
$G$-Banach algebra ({\it i.e.} a Banach algebra on which $G$ acts
continuously by isometric automorphisms $g : a\mapsto g(a)$). The space
$C_c(G,A)$ of $A$-valued continuous compactly supported functions on
$G$ is endowed with the following convolution product~:
$f*f'(g)=\int _G f(h)h(f'(h^{-1}g))dh$ and the completion
$L^1(G,A)$ of $C_c(G,A)$ for the norm $\|f\|=\int _G \|f(g)\|_A dg$ is
a Banach algebra.
More generally for any unconditional completion $\A(G)$, we define
$\A(G,A)$ to be the completion of $C_c(G,A)$ for the norm
$\|f\|_{\A(G,A)}=\big\|g\mapsto \|f(g)\|_A\big\|_{\A(G)}$.

For any $G$-Banach algebras $A$ and $B$,
we define in~\cite{these} an abelian group
$KK_G^\ban(A,B)$. This is a contravariant functor in $A$ and a covariant
functor in $B$. When $G=1$ this is equal to $KK^\ban(A,B)$. For any
unconditional completion $\A(G)$ of $C_c(G)$,  there is descent morphism
$KK_G^\ban(A,B)\to
KK^\ban(\A(G,A),\A(G,B))$.

These constructions are adaptations of the classical constructions for
$C^*$-algebras : for any $G$-$C^*$-algebra $A$
({\it i.e.} $G$ acts continuously by $C^*$-algebras automorphisms on $A$)
we have  a natural $C^*$-algebra $C^*_{\red}(G,A)$ containing $L^1(G,A)$ as
a dense subalgebra. If $B$ is another $G$-$C^*$-algebra,
 Kasparov defined an abelian group
$KK_G(A,B)$. This is a contravariant functor in $A$ and a covariant
functor in $B$. When $G=1$ this is equal to $KK(A,B)$.
There is an associative and distributive  product $KK_G(A,B)\otimes
KK_G(B,C)\to KK_G(A,C)$ and a descent morphism $KK_G(A,B)\to
KK(C^*_{\red}(G,A),C^*_{\red}(G,B))$.

Let $K^G_*(\underline{E}G,A)$, $*=0,1$, be the inductive limit over
$G$-invariant $G$-compact subsets $Z$ of $\underline EG$ of
$KK_{G,*}(C_0(Z), A)$. Then the assembly map
$$\mu_{\mathrm{red},A}:K^G_*(\underline{E}G,A)\to
K_*(C^*_{\mathrm{red}}(G,A))$$
 is defined in~\cite{baumconnes} and similar maps  $\mu_{L^1,A}$,
and more generally $\mu_{\A,A}$ for any unconditional completion
$\A(G)$, can be  defined.

The  Baum-Connes conjecture ``with coefficients'' claims that
$\mu_{\mathrm{red},A}$ is an isomorphism and the Bost conjecture ``with
coefficients'' claims that $\mu_{L^1,A}$ is an isomorphism.
Theorems~\ref{gammainj},~\ref{bostconj},~\ref{uncondthm} are still true with
arbitrary coefficients.

The surjectivity of the Baum-Connes conjecture with coefficients has
been counter-exampled recently (Higson, Lafforgue, Ozawa, Skandalis, Yu)
using a random group constructed by Gromov (\cite{spacequest})
but Bost conjecture with coefficients still stands.
If the Baum-Connes conjecture with coefficients is true for a group,
it is true also for all its closed subgroups; the Baum-Connes
conjecture with coefficients is also stable under various kinds of
extensions (Chabert~\cite{chabert},
Chabert-Echterhoff~\cite{chabertechterhoff},  Oyono~\cite{oyono}, and
Tu~\cite{tutrees}).

Kasparov's equivariant KK-theory was generalized to groupoids by Le
Gall \cite{kaspkk,legall,legall1}
and this generalized KK-theory was applied by Tu
in~\cite{tufeuilletages,tufhyp} to
the bijectivity  of the Baum-Connes map for amenable groupoids and the
injectivity for (the holonomy groupoids of)
hyperbolic foliations.
It is possible to generalize also Banach KK-theory and unconditional
completions. In this way we obtain the Baum-Connes conjecture  for any
hyperbolic group, with coefficients in any commutative $C^*$-algebra,  and also
for foliations with compact basis, admitting a (strictly) negatively curved
longitudinal riemannian metric, and such that the holonomy groupoid is
Hausdorff and has simply connected fibers (not yet published).


\section{Discrete series representations of connected \\ semi-simple Lie
groups}

\vskip-5mm \hspace{5mm}

In this part we examine how the Baum-Connes conjecture for a connected
semi-simple Lie group with finite center can be used to establish the
construction of the discrete series by Dirac induction
(\cite{harishchandra,partha,atiyahschmid}).  That this is morally true is
known from the beginning of the conjecture (see for
instance~\cite{kingston}). In the proof we shall introduce
3 ingredients : these are classical facts stated here without proof.
Parts of the argument apply to more general groups (not connected, not
semi-simple).

This work owes its existence to Paul Baum. He asked me to study the
problem and we discussed a lot.

\subsection*{2.1. Dirac operators} \addsec

\vskip-5mm \hspace{5mm}

Let $G$ be a Lie group, with a finite number of connected
components, and $K$ a maximal compact
subgroup. We assume that there exists a $G$-invariant orientation on
$G/K$. For the sake of simplicity, we assume that $G/K$ admits a
$G$-invariant spin structure (it is true anyway for a two fold
covering of $G$). More precisely let $\mathfrak p$ be a complementary
subspace for the Lie algebra $\mathfrak k$ of $K$ in the  Lie algebra
$\mathfrak g$ of $G$. We choose $\mathfrak p$ such that it is invariant for the
adjoint action of $K$ and we endow it with a $K$-invariant
euclidian metric. The above assumption means that the homomorphism $K\to
\mathrm{SO}(\mathfrak p)$ lifts to $\mathrm{Spin}(\mathfrak p)$. We denote by
$S$ the associated
spin representation of $K$. If $\dim (G/K)$ is even, $S$ is
$\Z/2\Z$-graded. We write $i=\dim (G/K)\ \ [2]$.

We denote by $R(K)$ the (complex) representation ring of $K$ and for any
finite dimensional representation $V$ of $K$ we denote by $[V]$ its
class in $R(K)$.

Let $V$ be a finite dimensional representation of $K$. Let $E_V$ be the
right Banach (in fact Hilbert) module over $C^*_{\mathrm{red}}(G)$
($\Z/2\Z$-graded if $i=0\ [2]$) whose
elements are the $K$-invariant elements in $V^*\otimes S^*\otimes
C^*_{\mathrm{red}}(G)$, where $K$ acts by left translations  on
$C^*_{\mathrm{red}}(G)$. Let
$D_V$ be the unbounded $C^*_{\mathrm{red}}(G)$-linear
operator on $E_V$ equal to $\sum 1\otimes c(p_i)\otimes p_i$, where
the sum is over $i$, $(p_i)$ is an orthonormal basis of $\mathfrak p$,
$p_i$ denotes also the associated right invariant vector field on $G$,
and $c(p_i)$ is the Clifford multiplication by $p_i$. Let
$T_V=\frac{D_V}{\sqrt{1+D_V^2}}$. Then we define $[d_V]\in
K_i(C^*_{\mathrm{red}}(G))$ to be the class of the Fredholm module
$(E_V,T_V)$ over $C^*_{\mathrm{red}}(G)$.

In other words, $E_V$ is the completion  of the space of
smooth compactly supported sections of the bundle on $K\backslash G$ associated
to the representation $V^*\otimes S^*$ of $K$, for the
norm $\|w\|=\sup _{f\in L^2(G),
\|f\|_{L^2(G)}=1}\|w*f\|_{L^2((V^*\otimes S^*)\times _K G)}$,  and $D_V$
is the Dirac operator, twisted by $V^*$.

\noindent{\bf Connes-Kasparov conjecture.} {\it The group morphism
$\mu_{\red}:R(K)\to K_i(C^*_{\mathrm{red}}(G))$ defined by
$[V]\mapsto [d_V]$ is an isomomorphism, and
$K_{i+1}(C^*_{\mathrm{red}}(G))=0$.}

This is a special case of the Baum-Connes conjecture because we may
take $\underline{E}G=G/K$ and thus $K_i^G(\underline{E}G)=R(K)$ and
$K_{i+1}^G(\underline{E}G)=0$.
It was checked for $G$ connected reductive linear in~\cite{wassermann}
and the Baum-Connes conjecture was proved for any reductive group
in~\cite{these} (see c) of the corollary~\ref{maincor} above).

The following lemma has been suggested to me by Francois Pierrot.
Assume that $i$ is even. Let
moreover $H$ be a unitary tempered admissible representation of
$G$. This implies that we have a $C^*$-homomorphism $C^*_{\mathrm{red}}(G)\to
\K(H)$. For any element $x\in K_0(C^*_{\mathrm{red}}(G))$ we denote by $\langle
H,x\rangle\in \Z$ the image of $x$ by $K_0(C^*_{\mathrm{red}}(G))\to
K_0(\K(H))=\Z$.  If $x$ is the class of an idempotent $p\in
C^*_{\red}(G)$, the image of $p$ in $\K(H)$ is a finite rank
projector, whose rank is $\langle H,x\rangle$.

\begin{lemma} We have $\langle H,[d_V]\rangle=\dim(V^*\otimes
S^*\otimes H)^K$.
\end{lemma}

\subsection*{2.2. Dual-Dirac operators} \addsec

\vskip-5mm \hspace{5mm}

From now on we assume that $G$ is a connected semi-simple Lie
group with finite center and we still assume that $G/K$ has a
$G$-invariant spin structure. Kasparov has constructed an element
$\eta \in \mathrm{Hom}(K_i(C^*_{\mathrm{red}}(G)), R(K))$ (coming from
an element of $KK_i(C^*_{\mathrm{red}}(G),C^*_{\mathrm{red}}(K))$,
itself coming from an element of $KK_{G,i}(\C,C_0(G/K))$). Kasparov
has shown that $\eta\circ
\mu_{\red}=\mathrm{Id}_{R(K)}$~\cite{conspectus,kaspkk}.

Here is the detail of the construction.
The $G$-invariant riemannian structure on
$G/K$ given by the chosen
$K$-invariant euclidian metric on $\mathfrak p$ has  non-positive
curvature. Let $\rho$ be the distance to the origin and $\xi=d(\sqrt
{1+\rho ^2})$. Let $V$ be a finite dimensional complex representation of
$K$, endowed
with an invariant hermitian metric. Let $H_V$ be the space of $L^2$
sections of the hermitian $G$-equivariant fibre bundle on $G/K$
associated to the representation of $K$ on $S\otimes V$ and
let $c_{\xi,V}$ be the Clifford multiplication by $\xi$. In other words $H_V$
is the subspace of $K$-invariant vectors in $L^2(G)\otimes S\otimes
V$, where $K$ acts  by right translations $L^2(G)$, and $c_{\xi,V}$ is
the restriction to this subspace of the tensor
product of the Clifford multiplication by $\xi$ on $L^2(G)\otimes S$
with $\Id_V$. Left translation by $G$ on $G/K$ or on $L^2(G)$ gives
rise to a ($C^*$-)morphism $\pi_V:C^*_{\mathrm{red}}(G)\to \L_\C(H_V)$ and
$(H_V,\pi_V,c_{\xi,V})$ defines $\eta_V\in
KK_i^\ban(C^*_{\mathrm{red}}(G),\C)$ (in
fact in $KK_i(C^*_{\mathrm{red}}(G),\C)$). We denote
by $[\eta_V]\in \mathrm{Hom}(K_i(C^*_{\mathrm{red}}(G)),\Z)$ the
associated map,
and $\eta =\sum _{V} [\eta_V][V]\in
\mathrm{Hom}(K_i(C^*_{\mathrm{red}}(G)),R(K))$, where the sum is over the
irreducible representations of $K$.

Since the Connes-Kasparov conjecture is true,  $\mu_{\red}:R(K)\to
K_i(C^*_{\mathrm{red}}(G))$
and $\eta:K_i(C^*_{\mathrm{red}}(G))\to R(K)$ are inverse of each other and
$K_{i+1}(C^*_{\mathrm{red}}(G))=0$.

Let $H$ be a discrete series representation of $G$, {\it i.e.} an
irreducible unitary representation with a positive mass in the Plancherel
measure. We recall that this is equivalent to the fact that some (whence
all) matrix coefficient $c_x(g)=\langle x,\pi(g) x\rangle$, $x\in H$,
$\|x\|=1$, is square-integrable. Then $\|c_x\|_{L^2(G)}^2$ is
independant of $x$, and its inverse is the formal degree $d_H$ of $H$, which
is also the mass of $H$ in the Plancherel measure.
We introduce a first ingredient.

{\bf Ingredient $1$.} All discrete series representations of $G$ are
isolated in the tempered dual.

In other words, all matrix coefficients belong to
$C^*_{\mathrm{red}}(G)$. In fact a
standard asymptotic expansion argument shows that for any $K$-finite
vector $x\in H$, $c_x$  belongs to the Schwartz algebra
(\cite{harishchandra}, II, corollary 1 page 77).

Therefore there  exists an
idempotent $p\in C^*_{\mathrm{red}}(G)$ such that the image in $L^2(G)$ of the
image of $p$ by the left regular representation is $H^*$ as a
representation of $G$ on the right. In fact we can take
$p=d_H\overline{c_x}$ for
any $x\in H$, $\|x\|=1$, where $\overline{c_x}(g)=\overline{c_x(g)}$. The
class of $p$ in $K_0(C^*_{\mathrm{red}}(G))$ only
depends on $H$ and we denote it by $[H]$. It is easy to see that
$i: \oplus_H \Z \to K_0(C^*_{\mathrm{red}}(G))$, $(n_H)_H\mapsto
\sum_H n_H [H]$, where the sums are
over the discrete series representations  of $G$, is an injection.
Indeed, if $H$ and $H'$ are discrete
series representations of $G$, $\langle H',[H]\rangle=1$ if $H=H'$ and
$0$ otherwise.

As a corollary we see that if $i=1 \ [2]$, $G$ has no discrete series
representations. From now on we assume $i=0\ [2]$.

The first part of the following lemma was suggested to me by Georges
Skandalis. Let $H$ be a discrete series representation of $G$. We
write $\eta ([H])=\sum_V n_V [V]$ in
$R(K)$ where the sum is finite and over the irreducible
representations of $K$ (in the notation above, $n_V=[\eta_V]([H])$).

\begin{lemma}
If $V$ is an irreducible representation of $K$,
$n_V=\dim(H^*\otimes S\otimes V)^K$ and therefore
$n_V=\langle H,[d_V]\rangle$.
\end{lemma}

We have
$1=\langle H,[H]\rangle =\langle H,\mu_{\red}\circ \eta ([H])\rangle=
\sum_V n_V\langle H,[d_V]\rangle=\sum_V n_V^2$.
Therefore one of the $n_V$ is $\pm 1$ and the others are $0$.

Alternatively we can consider the morphisms
\begin{eqnarray*}\oplus_V \Z [V]=R(K)\overset{\mu_\red}{\to}
K_0(C^*_{\mathrm{red}}(G))\overset{\pi}{\to}\prod_H \Z
\text{ where } \pi(x)=(\langle H,x\rangle )_H
\\ \text{
and }\ \ \ \ \ \
\oplus_H \Z \overset{i}{\to}
K_0(C^*_{\mathrm{red}}(G))\overset{\eta}{\to}R(K)=\oplus_V \Z
[V]\end{eqnarray*}
where the sums are over the irreducible representations $V$ of $K$ and
the discrete series representations $H$ of $G$.
 Their product  $\pi\circ \mu_\red\circ \eta \circ i=\pi\circ
i$ is equal to the inclusion of $\oplus_H \Z $ in $\prod_H \Z $
 and their matrices in the
base $([V])_V$ and the canonical base of $\oplus_H \Z$
  are transpose of each other. Therefore each
column of the matrix of $\eta\circ i$ contains exactly one non-zero
coefficient, which is equal to $\pm 1$.
A posteriori, $\pi$ takes its values in $\oplus_H \Z $.

\begin{corollary}
The discrete series representations of $G$ are in bijection with a
subset of the set of isomorphism classes of irreducible
representations of $K$. The
irreducible representation $V$ of $K$ associated to a discrete series
representation $H$ is such that $V=\pm (H\otimes S^*)$ as a formal
combination of irreducible representations of $K$, and $H$ occurs in
the kernel of the twisted Dirac operator $D_V$.
\end{corollary}

\begin{corollary} If $\mathrm{rank}\, G\neq \mathrm{rank}\, K$, $G$ has
no discrete series.
\end{corollary}

In this case $S^*$ is $0$ in $R(K)$ (Barbasch and
Moscovici~\cite{barbasch} (1.2.5) page 156) : this was indicated
to me by Henri Moscovici.

\subsection*{2.3. A trace formula} \addsec

\vskip-5mm \hspace{5mm}

From now on we assume that $\mathrm{rank}\, G= \mathrm{rank}\, K$. Let
$T$ a maximal torus in $K$ (therefore also in $G$). Choose a Weyl
chamber for the root system of $\mathfrak g$ and choose the Weyl chamber
of the root system of $\mathfrak k$ containing it. Let $V$ be an
irreducible representation of $K$, $\mu$ its highest wheight, and
$\lambda=\mu +\rho_K$ where $\rho_K$ is the half sum of the positive
roots of $\mathfrak k$.

We recall that the unbounded trace $\mathrm{Tr}:C^*_{\mathrm{red}}(G)\to \R,
f\mapsto f(1)$ gives
rise to a group morphism $K_0(C^*_{\mathrm{red}}(G))\to \R$. When $H$ is a discrete series
representation of $G$, $ \mathrm{Tr}([H])$ is the value at $1$ of
$p=d_H\overline{c_x}$ for some $x\in H$, $\|x\|=1$, and therefore it is the
formal degree $d_H$ of $H$ and is $>0$.

{\bf Ingredient 2.} $\mathrm{Tr}([d_V])=\prod _{\alpha \in
\Psi}\frac{(\lambda ,\alpha)}{(\rho,\alpha)}$,
where $\Psi$ is the set of simple roots of the chosen positive root
system in $\mathfrak g$, and $\rho$ is the half sum of the positive roots
of this system.

In this formula is used a right normalization of the Haar measure (if
$G$ is linear it is the one for which the maximal compact subgroup of
the complexification of $G$ has measure $1$).
This formula is proven in~\cite{connesmoscovici} by a heat equation
method, and in~\cite{atiyahschmid} by Atiyah's $L^2$-index theorem.

\begin{corollary}
If $\lambda$ is singular for $\mathfrak g$, $[V]$ does not correspond to a
discrete series representation of $G$.
\end{corollary}

{\bf Ingredient 3.} For any $x\in K_0(C^*_{\mathrm{red}}(G))$ such that
$\mathrm{Tr}(x)\neq 0$, there is a discrete series representation $H$
such that $\langle H,x\rangle\neq 0$.

By the Plancherel formula, if $\hat G$ is the tempered spectrum of
$G$, $\mathrm{Tr}(x)=\int_{\hat G}\langle H,x\rangle\ dH$. We have to
prove that, for almost all $H$ outside the discrete series, $\langle
H,x\rangle=0$. There are several possible arguments :
\begin{itemize}
\item almost all $H$ outside the discrete series are induced
from a parabolic subgroup and belong to a family of representations
indexed by some $\R^p$, but $\langle H',x\rangle$ is constant when
$H'$ varies in this family and goes to $0$ when $H'$ goes to
infinity,
\item write $x=[d_V]$ for some $V$, then the $H$ outside the
discrete series with  $\langle
H,x\rangle\neq 0$ have measure $0$  by~\cite{atiyahschmid} p15
(3.19), p50 (9.8) and p51 (9.12) or by~\cite{connesmoscovici}
p318-320.
\end{itemize}

\begin{corollary}If $\lambda$ is not singular for $\mathfrak g$, $[V]$ does
correspond to a discrete series representation, whose formal degree is
$\Big| \prod _{\alpha \in
\Psi}\frac{(\lambda ,\alpha)}{(\rho,\alpha)}\Big|$.
\end{corollary}

We have recovered some results proved in~\cite{harishchandra},
\cite{partha}  and~\cite{atiyahschmid}.

\label{lastpage}

\end{document}